# On Parameterization of Some Problems in Number Theory

Sergey N. Artekha
Space Research Institute of RAS, 117997 Moscow, Russia
E-mail: sergey.arteha@gmail.com

**Abstract**

The article is devoted to the parametric solution of some Diophantine equations with two variables. In essence, we find an additional parametric relationship between the solution of the original equation and the simpler linear equation $ax-by = \pm 1$. The method of hidden parameters helps to obtain a parametric representation of integer solutions for a number of problems in number theory. The method consists in that, first, instead of each quantity included in the original equation, we introduce an additional number of parameters using a linear or non-linear function, which allows us to find both connections between the parameters and some parametric solution of the desired equation. Then we reset (nullify) the introduced linear or non-linear functional dependence (hidden parameters), which in the original formulation of the problem would mean just changing the letter designation of the parameters. However, as a result of these actions, a simplified parametric solution is obtained and the found relationship between the parameters remains. In this paper, the method is discussed in detail using the example of parametric solving such Diophantine equation as the Pell equation. The method of such parameterization is also applied to problems related to Pell's equation. The methods of parameterization are discussed for the Delaunay equation, the Nagell equation and similar equations of higher powers.

**Keywords:** Diophantine equations, Pell equation, Delaunay equation, Nagell equation, parametric solution

## 1 Introduction

In the article [1], a method of hidden parameters was proposed for solving the Pell equation in natural numbers. To significantly reduce the time of numerical calculations for specific numbers of $A$, a number of simplifications were carried out in this work for five different special cases, into which the general problem can be divided. Also, for these cases, algorithms for quickly filling out the solution table for the sequence of numbers $A$ have been developed. In this paper, we are not interested in reducing the calculation algorithm for fixed particular cases. Let us recall the essence of the method for its further application in order to obtain parametric representations for the solutions of some problems of number theory. Instead of each quantity included in the original equation, we introduce an additional number of parameters using a linear or non-linear function, which allows us to find both relationships between the parameters and some parametric solution of the desired equation. After that, we can reset (nullify) the introduced additional linear or non-linear dependence (in fact, hidden parameters), which in the original formulation of the problem would mean just replacing the letter designation of the parameters. However, as a result of these actions, a simplified parametric solution is obtained and the found relationship between the parameters remains. But we can remain (without nullify) the introduced additional parameters, then the solution turns out to be more cumbersome, with a large number of parameters, but capturing a larger number of parametric solutions. As a result, we find a parametric relationship between the solution of the original equation and the simpler equation $ax-by = \pm 1$. We will demonstrate in detail a method for solving the Pell equation in a way different from that used in [1].

The main goal of this article is to develop the proposed method in the most general form, including applying it to higher-degree Diophantine equations [2-10], such as the Delaunay equation, the Nagell equation and similar equations, and to obtain parametric representations for their solutions.

## 2 Pell's equation and related topics

We consider another direct representation of the solution to Pell's equation

$$Ax^2 + 1 = y^2 \tag{1}$$

for the case

$$A = p_{01}p_{02}\,,\ x = 2fg, \tag{2}$$

where all the variables introduced for factoring denote some coprime numbers (but they can turn out to be anything: prime, composite, and one). We rewrite (1) in the form of a system

$$\begin{aligned} y + 1 = 2p_{01}f^2, \\ y - 1 = 2p_{02}g^2, \end{aligned} \tag{3}$$

from which you can take any of the two expressions for the value of $y$ and write the following condition:

$$p_{01}f^2 - p_{02}g^2 = 1. \tag{4}$$

If for fixed values of $f$ and g (that is, for a fixed value of $x$) we find one solution to (4), then we will know all values of $A$ (there are an infinite number of them), for which the same value of $x$ will also be a solution:

$$A^i = p_1^i p_2^i\,,\ \ p_1^i = p_{01} + ig^2\,,\ \ p_2^i = p_{02} + if^2\,,\ \ i \in Z. \tag{5}$$

We will make linear substitutions in equation (4):

$$p_{01} = ak + b,\ p_{02} = ck + d,\ f = rk + s,\ g = tk + u. \tag{6}$$

We group the resulting equation according to degrees of $k$, and so that it is performed with any values of $k$, equate each expression with the appropriate degree of $k$ to zero. To resolve the resulting system of three equations with senior degrees in integers, we make substitutions $t = zt_0$ , $r = zr_0$ and get a solution:

$$a = 2zr_0t_0^3\,,\ b = t_0^2(3ur_0 - st_0)\,,\ c = 2zr_0^3t_0\,,\ d = r_0^2(3st_0 - ur_0), \tag{7}$$

and from the equation of the system with the zero degree of $k$, the condition follows:

$$ur_0 - st_0 = 1. \tag{8}$$

Thus, we received the connection of equation (1) for case (2) with the simplest linear equation (8). Note that the obtained solution is true for any $k$, including $k = 0$. Although in this particular case, our substitutions (6) would only mean simple changing the letter designation of variables (hidden parameters).

If, with fixed values of $u$ and $s$, we find the solution $r_0$ and $t_0$ to the equation (8), then the solutions will also be

$$r_j = r_0 + js\,,\;\; t_j = t_0 + ju\,,\;\; j \in Z. \qquad (9)$$

Then the final relationship between the solution of the equation (8) and the solution (2) of the equation (1) is as follows:

$$x = 2\big((r_0 + js)k + s\big)\big((t_0 + ju)k + u\big), \qquad (10)$$

$$A = \Big(2k(r_0 + js)(t_0 + ju)^3 + (t_0 + ju)^2(3u(r_0 + js) - s(t_0 + ju)) + i\big((t_0 + ju)k + u\big)^2\Big) \times$$
$$\times \Big(2k(r_0 + js)^3(t_0 + ju) + (r_0 + js)^2(3s(t_0 + ju) - u(r_0 + js)) + i\big((r_0 + js)k + s\big)^2\Big). \qquad (11)$$

Here *i*, *j*, *k* are arbitrary integers. However, note that there are also some rational numbers that lead to natural values *x*, *A*. (This can be done, if we are interested in the question, for which values of *A* some fixed value of *x* is a solution.)

Thus, we obtain the following

**Proposition 1** *All solutions to equation (1) for the case (2) are given by expressions (10), (11), where the parameters u and s are related by relation (8), and i, j, k are arbitrary integers.*

We can write another solution if, on the contrary, we fix $r_0$ and $t_0$ and find solutions $u_0$ and $s_0$ of the equation (8); the solutions also will:

$$u_j = u_0 + jt_0\,,\;\; s_j = s_0 + jr_0\,,\;\; j \in Z. \qquad (12)$$

As a result, another ultimate relationship between the solution of equation (8) and the solution (2) of equation (1) looks as follows:

$$x = 2(r_0k + s_0 + jr_0)(t_0k + u_0 + jt_0), \qquad (13)$$

$$A = (2kr_0t_0^3 + t_0^2(3(u_0 + jt_0)r_0 - (s_0 + jr_0)t_0) + i(t_0k + u_0 + jt_0)^2) \times$$
$$\times (2kr_0^3t_0 + r_0^2(3(s_0 + jr_0)t_0 - (u_0 + jt_0)r_0) + i(r_0k + s_0 + jr_0)^2). \qquad (14)$$

Here *i*, *j*, *k* are again arbitrary integers (but there are also some rational numbers leading to integers *x*, *A*).

Thus, we obtain the following

**Proposition 2** *All solutions to equation (1) for the case (2) can be presented by expressions (13), (14), where the parameters* $u_0$ and $s_0$ *are related by relation (8), and i, j, k are arbitrary integers.*

As a solution to equation (8), you can choose any identity, for example:

$$n_1 \cdot n_1n_2^2 - (n_1n_2 - 1) \cdot (n_1n_2 + 1) \equiv 1, \qquad (15)$$

or more general identities, for example:

$$\left[\begin{array}{l} n_1^k - (n_1 - 1)\big(n_1^{k-1} + \cdots + 1\big) \equiv 1 \\ (n_2 + 1)\big(n_2^{2k} - \cdots + 1\big) - n_2^{2k+1} \equiv 1. \end{array}\right. \qquad (15^*)$$

We choose, which of the multipliers from the identity (15) (or from any similar identity, for example, (15*)) to assign for which of variables in (8). After that, we substitute them into the solutions (10), (11) or (13), (14). As a result, we obtain parametric solutions of Pell's equations. When substituting a set of specific numerical parameters into the resulting identic formula, some solutions may be repeated, and for a particular *A*, several different solutions may appear (of course, not all of these solutions turn out to be primitive - with a minimal non-zero natural *x*). Sequential filling of the solution table occurs quite fast.

A similar approach can be used to obtain parametric solutions not only for even solutions *x* from (2), but also for odd $x = fg$, which together covers all possible solutions. In this particular case, in equation (4), the number 2 will be on the right, and all solutions in (7) must be doubled. Then equation (8) turns out to be unchanged. All substitutions retain their former form, and the corresponding analogous solutions will be:

$$x = \big((r_0 + js)k + s\big)\big((t_0 + ju)k + u\big), \tag{10*}$$

$$A = \Big[4k(r_0 + js)(t_0 + ju)^3 + 2(t_0 + ju)^2\big(3u(r_0 + js) - s(t_0 + ju)\big) + \\ + i\big((t_0 + ju)k + u\big)^2\Big]\,[4k(r_0 + js)^3(t_0 + ju) + 2(r_0 + js)^2 \\ (3s(t_0 + ju) - u(r_0 + js)) + +i\big((r_0 + js)k + s\big)^2\Big]; \tag{11*}$$

$$x = (r_0k + s_0 + jr_0)(t_0k + u_0 + jt_0), \tag{13*}$$

$$A = (4kr_0t_0^3 + 2t_0^2(3(u_0 + jt_0)r_0 - (s_0 + jr_0)t_0) + i(t_0k + u_0 + jt_0)^2) \times \\ \times (4kr_0^3t_0 + 2r_0^2(3(s_0 + jr_0)t_0 - (u_0 + jt_0)r_0) + i(r_0k + s_0 + jr_0)^2). \tag{14*}$$

The condition will coincide with equation (8) and all subsequent substitutions, for example, from the identity (15), (15*) or similar identities, remain fully valid.

Thus, we write the following

**Proposition 3** *All solutions to equation (1) for the case of* odd $x = fg$ *can be presented by expressions (10*), (11*) and (13*), (14*), where the parameters* $u_0$, $s_0$, $r_0$ and $t_0$ *are related by relation* $u_0r_0 - s_0t_0 = 1$*, and i, j, k are arbitrary integers.*

Similarly, with the help of linear substitutions, one can search for a solution to the equation

$$ax^2 - by^2 + cxy = n$$

with a given right side and a fixed one of the coefficients on the left side.

Let us now turn to higher-degree equations. If the equation can be reduced to an equation of type (4) with linear coefficients, then its further solution is elementary. Let, for example, the same power *n* be present in the equation:

$$p_{01}f^n - p_{02}g^n = 1. \tag{16}$$

For example, the solution of the equation

$$Ax^n + 1 = y^2 \tag{17}$$

for the following two cases:

$$x = 2fg,\quad \begin{bmatrix} 1)\ A = p_{01}^*p_{02}\,,\ p_{01} = 2^{n-2}p_{01}^* \\ 2)\ A = p_{01}p_{02}^*\,,\ p_{02} = 2^{n-2}p_{02}^* \end{bmatrix} \tag{18}$$

reduces to such the equation (16), where all the variables introduced for factoring denote some coprime numbers. By analogy with quadratic equations, it is also easy to reduce the search for odd solutions *x* to the solution of the equation (16).

If for fixed values of *f* and *g* we find one (any) solution to (16), then we will know all solutions (all values of these quantities):

$$p_1^i = p_{01} + ig^n\,,\ \ p_2^i = p_{02} + if^n\,,\ \ i \in Z. \tag{19}$$

To find a solution, the following substitutions can be made in (16):

$$p_{01} = b + \sum_{j=1}^{n-1} a_jk^j\,,\ p_{02} = d + \sum_{j=1}^{n-1} c_jk^j\,,\ f = rk + s,\ g = tk + u. \tag{20}$$

Next, you need to group the resulting equation by powers of $k$, and in order for it to be fulfilled for any values of $k$, equate each expression before the corresponding power of $k$ to zero. After that, the system of linear equations is solved relative all coefficients in $p_{01}$ and $p_{02}$.

However, in order to find one solution (after all, this is enough), you can do it easier. If we want to find a connection between the solutions of equation (16) and equation

$$uf - sg = 1, \tag{21}$$

then we raise the last equation to the power $2n - 1$ and obtain the required quantities:

$$\begin{cases} p_{01} = u^{2n-1}f^{n-1} - C_{2n-1}^{1}u^{2n-2}f^{n-2}sg + \cdots + (-1)^{n-1}C_{2n-1}^{n-1}u^{n}s^{n-1}g^{n-1}, \\ p_{02} = s^{2n-1}g^{n-1} - C_{2n-1}^{1}ufs^{2n-2}g^{n-2} + \cdots + (-1)^{n-1}C_{2n-1}^{n-1}u^{n-1}f^{n-1}s^{n}, \end{cases} \tag{22}$$

written in terms of the expansion coefficients of the Newton binomial, taking into account the symmetry of the coefficients. For the first case in (18), the number $u$ must be even. For the second case in (18), the number $s$ must be even. In any of these cases the number $i$ in (19) must be $i = 2^{n-2}i_1$. The choice of some identities, for example, of the type (15) or (15*) as solutions to the equation (21), allows us to write the solutions of the equations (16) in a parametric form. Obviously, we can use substitutions that do not change equation (21), for example: $u \to u + jg, s \to s + jf$, or $f \to f + ks, g \to g + ku$, where $j$, $k$ are arbitrary integers.

Thus, we obtain the following

**Proposition 4** *All solutions to equations (16) and (17) for the case (18) are given by expressions (19), (22), where the parameters u and s are related by relation (21), and i is an arbitrary integer.*

In essence, we have obtained all solutions to the following problem in integers: for what numbers $A$, a given even $x$ is a solution to equation (17).

Similarly, you can act in the case of an equation with different powers:

$$p_{01}f^{n} - p_{02}g^{k} = 1. \tag{23}$$

For fixed values of $f$ and $g$ (knowing one solution), all solutions have the form:

$$p_1^i = p_{01} + ig^k\ ,\ \ p_2^i = p_{02} + if^n\ ,\ \ i \in Z, \tag{24}$$

and the relationship with the solutions of equation (21) after raising to the power n + k − 1 is given by the expressions:

$$\begin{cases} p_{01} = u^{n+k-1}f^{k-1} - C_{n+k-1}^{1}u^{n+k-2}f^{k-2}sg + \cdots + (-1)^{k-1}C_{n+k-1}^{k-1}u^{n}s^{k-1}g^{k-1}, \\ p_{02} = -(-s)^{n+k-1}g^{n-1} + (-s)^{n+k-2}C_{n+k-1}^{1}ufg^{n-2} - \cdots - (-s)^{k}C_{n+k-1}^{k}u^{n-1}f^{n-1}, \end{cases} \tag{25}$$

which are written in terms of the expansion coefficients of the Newton binomial. Again, choosing as solutions to the equation (21) some identities, for example, of the type (15) or (15*), we can write the solutions of the equations (23) in a parametric form.

Therefore, we get the following

**Proposition 5** *All solutions to equation (23) are given by expressions (24), (25), where the parameters u and s are related by relation (21), and i is an arbitrary integer.*

Thus, if a pair of coefficients enters the original equation linearly (or the original equation can be reduced to a similar equation), then it is easy to find a connection between this equation and a linear equation of type (21) and obtain parametric solutions in symbolic form (in principle, cover all solutions). Note that in article [8], in fact, only some special cases of polynomial parametrization of the Pell equation are presented. If there is only one such a coefficient or the coefficients enter in a non-linear way, then the problem may become more complicated.

### 3 The Delaunay equation

We write the Delaunay equation in integers:

$$x^3 + By^3 = 1. \tag{26}$$

We consider the case when $y$ and $B$ are not divisible by 3 and write down an arbitrary partition into coprime factors (in the general case all values are integers):

$$B = p_{01}p_{02}\,,\ y = -fg. \tag{27}$$

As a result, we have a system of equations where the left-hand sides of the equations are coprime (the case when they have a common divisor of 3 can be considered in a completely similar way):

$$\begin{cases} x - 1 = p_{01}f^3 \\ x^2 + x + 1 = p_{02}g^3. \end{cases} \tag{28}$$

Multiplying the second equation by 4, we rewrite it as:

$$(2x+1)^2 + 3 \cdot 1^2 = 4p_{02}g^3. \tag{29}$$

Since there exists a single-valued division procedure for any form of the type

$$X^2 + 3 \cdot Y^2,$$

then all divisors on the right side of (29) must be represented in the same form:

$$p_{02} = a^2 + 3b^2, \quad g = q^2 + 3s^2. \tag{30}$$

Substituting these expressions and the expression $x$ from the first equation of the system (28) into the equation (29), we obtain the following system of equations:

$$\begin{cases} \pm\{2p_{01}f^3 + 3\} = A_0q[q^2 - 9s^2] \mp 9B_0s[q^2 - s^2], \\ \pm 1 = B_0q[q^2 - 9s^2] \pm 3A_0s[q^2 - s^2], \end{cases} \tag{31}$$

where the signs on the right-hand sides of both equations of the system (31) are consistent with each other, the signs ± on the left-hand sides of both equations are completely arbitrary (each is independent of anything), and the notations with consistent signs are also introduced:

$$A_0 = a \mp 3b, \ \ B_0 = a \pm b. \tag{32}$$

Let us include the sign ∓ from the right parts of the system (31) into the value s, then in the first equation of (31), the sign (+) will present on the right, and in the second equation of (31), the sign (− ) will present on the right. Let us also include the sign ± on the left-hand side of the first equation of (31) in the definition of the value *f*. As a result, we rewrite the system (31) in the form:

$$\begin{cases} A_0q[q^2 - 9s^2] + 9B_0s[q^2 - s^2] - 2p_{01}f^3 = \pm 3, \\ B_0q[q^2 - 9s^2] - 3A_0s[q^2 - s^2] = \pm 1. \end{cases} \tag{31*}$$

If one arbitrary solution $A_0$ and $B_0$ is known from the second equation of (31*) for the given values of $q$ and $s$, then we can write down all such solutions:

$$A_i = A_0 + i(q^3 - 9qs^2), \ \ B_i = B_0 + i(3q^2s - 3s^3), \ i \in Z, \tag{33}$$

and from (32) the values $a$ and $b$ can be expressed, i.e. all possible $p_{02}^i$ values can be found.

We are looking for a solution to the second equation of the system (31*) in the form:

$$q = hk + q_0, \ s = jk + s_0, \ A_0 = gk^2 + rk + t, \ B_0 = mk^2 + ck + d, \ k \in N. \tag{34}$$

Let us group the resulting equation by powers of $k$, and in order for it to be fulfilled for any values of $k$, we equate each expression before the corresponding power of $k$ to zero. We solve the resulting system of six equations in integers with respect to the variables *g*, *r*, *t*, *m*, *c*, *d* and find the quantities $A_0$ and $B_0$ of interest to us. In order for them to be integers, we equate the literal expression in the denominator of these quantities to ±1 and take it as some desired condition. In this case, to cancel 3 in the denominator of $A_0$, it is necessary to put $h = 3h_0$. Any other choice would lead to a contradiction with the desired condition when compared modulo 3. Taking into account this substitution, we write the final desired condition as follows:

$$jq_0 - 3h_0 s_0 = \pm 1. \quad (35)$$

To resolve this equation as a parametric identity, one can choose, similarly to (15) or (15*), any identity with tree terms and, comparing with (35), choose from them the values for the variables $j, q_0, h_0, s_0$. We also get the following result for the values of $A_0$ and $B_0$:

$$A_0 = \frac{81h_0^3k^2}{8}\left(9h_0^4 - 3h_0^2j^2 - 5j^4\right) +$$
$$+\frac{27k}{8}\left(18h_0q_0 + 15h_0^4j^2q_0 - j^6q_0 - 63h_0^5js_0 - 30h_0^3j^3s_0 - 3h_0j^5s_0\right) +$$
$$+\frac{9}{8}\left(9h_0^5q_0^2 + 10h_0^3j^2q_0^2 + 5h_0j^4q_0^2 - 15h_0^4jq_0s_0 - 30h_0^2j^3q_0s_0 - 3j^5q_0s_0 - -72h_0^5s_0^2\right), \quad (36)$$

$$B_0 = \frac{9jk^2}{8}\left(81h_0^6 + 45h_0^4j^2 + 3h_0^2j^4 - j^6\right) +$$
$$+\frac{9k}{8}\left(27h_0^5jq_0 + 30h_0^3j^3q_0 + 7h_0j^5q_0 + 81h_0^6s_0 - 15h_0^2j^4s_0 - 2j^6s_0\right) +$$
$$+\frac{1}{8}\left(8j^5q_0^2 + 243h_0^5s_0q_0 + 270h_0^3j^2s_0q_0 + 15h_0j^4s_0q_0 - 405h_0^4js_0^2 - -90h_0^2j^3s_0^2 - 9j^5s_0^2\right). \quad (37)$$

For all parameters, one can check all possible remainders modulo 8 and determine those linear sequences, with which expressions (36) and (37) become integers (the number 8 in the denominator is canceled). We will not write out all the sequences, since there are quite a lot of them. Note that in the method of hidden parameters, we can well restrict ourselves to the case $k = 0,\ q = q_0,\ s = s_0$, because it was enough for us to find at least one solution.

Let us substitute, instead of quantities $A_0$ and $B_0$, an arbitrary solution (33) into the first equation of the system (31*):

$$i(q^2 + 3s^2)^3 - 2p_{01}f^3 = A_0q[9s^2 - q^2] + 9B_0s[s^2 - q^2] \pm 3. \quad (38)$$

In this equation, the right-hand side is fixed, since *q* and *s* are given, and the values $A_0$ and $B_0$ are found in (36) and (37). On the left-hand side of the equation (38), the values *f*, *q*, and *s* will be given. The general solution (38) can be written as follows:

$$\begin{cases} i^u = I_0(A_0q[9s^2 - q^2] + 9B_0s[s^2 - q^2] \pm 3) + 2uf^3,\ u \in Z, \\ p_{01}^u = \frac{P_0(A_0q[9s^2-q^2]+9B_0s[s^2-q^2]\pm 3)}{2} + u(q^2 + 3s^2)^3, \end{cases} \quad (39)$$

where the quantities $I_0$ and $P_0$ (again, it is enough for us to find one solution) are found from the condition:

$$I_0(q^2 + 3s^2)^3 - P_0f^3 = 1. \quad (40)$$

A solution of this equation

$$\begin{cases} P_0 = w^3(w^2f^2 - 5wfQv + 10Q^2v^2), \\ I_0 = v^3(Q^2v^2 - 5Qvwf + 10w^2f^2), \end{cases} \quad (41)$$

with the introduced notation

$$Q = q^2 + 3s^2, \tag{42}$$

can be written through the solution of a more simple equation:

$$v(q^2 + 3s^2) - wf = 1. \tag{43}$$

To present this equation as a parametric identity, we can again choose some identities (see, for example, (15) or (15*), etc.) and select from them values for the variables $f, v, w$.

Thus, we finally obtain the following technique for writing nontrivial parametric solutions for the Delaunay equation:

**Algorithm 1**

(1) By specifying $q_0$ and $s_0$, we find the values of $j$ and $h_0$ from the equation (35), and
(2) from equalities (36) and (37) we obtain the values of $A_0$ and $B_0$.
(3) With the help of substitutions (33), we find all the values of $A_i$ and $B_i$, and
(4) from the system (32) all possible $a_i = (3B_i + A_i)/4$ and $b_i = \pm(B_i - A_i)/4$ follow, therefore,
(5) we have obtained all possible values of $p_{02}^i$ and $g^k$ from (30).
(6) Further, by setting the value of $f$, we find the values of $v$ and $w$ from the equation (43).
(7) The first expression in (41) gives the value of $P_0$, and
(8) from the second expression (39) we find all possible values of $p_{01}^u$.
(9) As the result, for a given value of $y$, we found all possible values of $B$ from (27). The corresponding $x$ values are found from the first equation (28).

Note that knowing a solution of the equation (35) allows us to find another, infinite sequence of non-trivial solutions:

$$\begin{cases} q_n = q_0 + 3nh_0, \\ s_n = s_0 + nj, \ n \in Z, \end{cases} \tag{44}$$

and then we proceed according to the method described in the previous paragraph.

In both approaches, filling the table of nontrivial solutions for different values of $B$ occurs elementarily. Note that for the parametric representation of solutions to the Delaunay equation, it is easy to search for solutions (35) in the form of identities, but it may be difficult to simultaneously find identical representations for solutions to equation (43).

### 4 Nagell's equation

We write the Nagell equation in integers:

$$x^4 - By^4 = \pm 1. \tag{45}$$

Here we consider the case with the upper (+) sign and write down an arbitrary partition into coprime factors:

$$B = P_{01}P_{02}\,, \ y = fg. \tag{46}$$

For the case of all odd factors, as a result, we have a system of equations:

$$\begin{cases} x^2 - 1 = P_{01}f^4, \\ x^2 + 1 = P_{02}g^4. \end{cases} \tag{47}$$

Since there exists a single-valued division procedure for a shape of the form

$$X^2 + Y^2,$$

then all divisors on the right-hand side of the second equation of the system (47) should be represented in such the form:

$$P_{02} = a^2 + b^2, \quad g = q^2 + s^2. \tag{48}$$

Substituting these expressions into the second equation of the system (47), we obtain the following solution *x* and the condition:

$$\begin{cases} x = 4aqs(q^2 - s^2) \pm b[(q^2 - s^2)^2 - (2qs)^2], \\ a[(q^2 - s^2)^2 - (2qs)^2] \mp 4bqs(q^2 - s^2) = \pm 1. \end{cases} \tag{49}$$

Here, the signs on the left-hand side of the condition are consistent with the signs in the solution, and the signs on the right-hand side of the condition are completely arbitrary (each does not depend on anything).

For the given values of *q* and *s*, if one arbitrary solution ($a_0$,$b_0$) is known from the second equation (49), then we can write down all such solutions:

$$a_i = a_0 \pm 4iqs(q^2 - s^2), \;\; \mathrm{b}_i = b_0 + i[(q^2 - s^2)^2 - (2qs)^2], \; i \in Z, \tag{50}$$

and from the first equation in (49) we obtain all possible $x_i$, and from (48) we can find all possible values of $P_{02}^i$. The second equation in (49) can be solved by the method of hidden parameters: by setting the values *q* and *s* as linear dependencies in *k* and looking for the values *a* and *b* in the form of quadratic dependencies; and then we separately equate the resulting expressions before different powers of *k* to zero (after finding the solution, you can choose *k* = 0). The result is the following solution:

$$a_0 = \frac{1}{4} qrs^2(7p^6 + 35p^4r^2 + 21p^2r^4 + 5r^6) -$$

$$-\tfrac{1}{4} pq^2s(5p^6 + 21p^4r^2 + 35p^2r^4 + 7r^6) + p^7s^3 - q^3r^7, \tag{51}$$

$$\begin{aligned} b_0 = \frac{\mp 1}{16} [pq^3(5p^6 + 21p^4r^2 + 35p^2r^4 + 35r^6) - q^2rs \times \\ (7p^6 + 35p^4r^2 + 105p^2r^4 + 29r^6) - pqs^2 \times \\ (29p^6 + 105p^4r^2 + 35p^2r^4 + 7r^6) +] \\ +rs^3(35p^6 + 35p^4r^2 + 21p^2r^4 + 5r^6). \end{aligned} \tag{52}$$

In the last expression, the $\mp$ signs are consistent with the first expression in (49) and expression (50). The values of $p$ and *r* are found from the following condition:

$$ps - qr = \pm 1, \tag{53}$$

in which the signs on the right-hand side are completely arbitrary. Again, all variables from this equation can be selected by comparison with parametric identities of the type (15), (15*), or others.

We introduce the following notations:

$$h = 4qs(q^2 - s^2), d = (q^2 - s^2)^2 - (2qs)^2, m = h^2 \pm d^2, N = a_0h \pm b_0d. \tag{54}$$

As a result, we write the solution in the form:

$$x_i = N + im, \quad P_{02}^i = a_i^2 + b_i^2. \tag{55}$$

However, we still need to coordinate this solution with the first equation of the system (47), from which two equations can be obtained:

$$\begin{cases} x - 1 = P_{11}f_1^4, \quad f = f_1f_2, \\ x + 1 = P_{12}f_2^4, \quad P_{01} = P_{11}P_{12}. \end{cases} \tag{56}$$

Subtracting the first equation from the second equation of (56), we obtain:

$$P_{12}f_2^4 - P_{11}f_1^4 = 2. \tag{57}$$

The solution to this equation in the general form will be:

$$\begin{cases} P_{12}^j = 2P_{12}^0 + jf_1^4, \\ P_{11}^j = 2P_{11}^0 + jf_2^4, \end{cases} \tag{58}$$

where $j \in Z$ and

$$\begin{cases} P_{12}^0 = 35c^4f_1^3r^3 - 21c^5f_1^2f_2r^2 + 7c^6f_1f_2^2r - c^7f_2^3, \\ P_{11}^0 = 35c^3f_2^3r^4 - 21c^2f_1f_2^2r^5 + 7cf_1^2f_2r^6 - r^7f_1^3, \end{cases} \tag{59}$$

and the values of $c$ and $r$ are found from the condition:

$$rf_1 - cf_2 = 1. \tag{60}$$

These variables can again be taken from a comparison with an arbitrary parametric identity of the type (15), (15*). Next, it follows from the first equation of (56) that

$$x_j = 1 + 2P_{11}^0f_1^4 + jf_1^4f_2^4. \tag{61}$$

Now we need to coordinate the solutions (61) and (55), i.e. find such values $i$ and $j$, for which the same solutions $x$ are obtained. Introducing the notations

$$E = 1 + 2P_{11}^0f_1^4 - N, \quad F = f_1^4f_2^4 \tag{62}$$

and equating (61) and (55), we obtain the equation

$$im - jF = E, \tag{63}$$

the solution of which in the general form (for $E \neq 0$) will be the system:

$$\begin{cases} i = Ei_0 + nF, \\ j = Ej_0 + nM, \end{cases} \tag{64}$$

$n \in Z$, and values of $i_0$ and $j_0$ can be found from the condition:

$$i_0m - j_0F = 1. \tag{65}$$

If it turns out that $E = 0$, then in this particular case, it will be $i = F$, $j = m$.

Thus, for writing parametric chains of solutions for the Nagell equation, we finally formulate

**Algorithm 2**

I. We give arbitrary values of $q$ and $s$ and find corresponding values of $p$ and $r$ from the equation (53) (even identical substitutions can be used at this stage).

II. After that, from the expressions (51) and (52), we obtain the values of $a_0$ and $b_0$, and from the definitions (54) we calculate the values of $h$, $d$, $m$ и $N$.

III. By setting arbitrary values of $f_1$ and $f_2$, we find the quantities $r$ and $c$ from the equation (60), and from the expressions (59) we find the quantities $P_{11}^0$ and $P_{12}^0$.

IV. Using the definitions (62) and solving the equation (65), we find values of $i_0$ and $j_0$.

V. From the expression (64), we find values of $i$ and $j$, and from the expression (55) we obtain values of $x_i$ and $P_{02}^i$.

VI. From the expression (58), we find values of $P_{11}^j$ и $P_{12}^j$, therefore, we know $P_{01}^j$.

VII. Using the remaining expressions from (56), (48) and (46), we find $B_i$ and $y_i$.

We will not consider cases where any factors in (46) are even - they can easily be resolved in a completely similar way (reduced to dividing some $P$ and a number of expressions obtained by some powers of 2).

If in the equation (45) with the minus sign on the right-hand side, instead of the fourth power of $x$, there would a second power: $x^2 - By^4 = -1$, then such an equation was already resolved using expressions (49) - (53). Moreover, the proposed technique allows one to analogously solve a similar equation with an arbitrary degree of $y$: $x^2 - By^n = -1$. But, unfortunately, the case with the minus sign on the right-hand side of the original equation (45) leads to the need to look for squares in some linear sequence or look for the possibility of its integer division by the fourth power of a certain number, which greatly complicates the search for a general parametric representation of the solution. A similar difficulty with finding the possibility for a linear sequence to be divided by the fifth power of a certain number has the solving of an equation of the fifth degree in $x$ and $y$, analogous to the Delaunay and Nagell equations: $x^5 - By^5 = -1$.

**5 Conclusion**

Thus, a method has been proposed that makes it possible to obtain various parametric representations of solutions to Diophantine equations. Using factorization and/or quadratic representations, we break the original equation into equations that have a pair of unknown free variables. And then we look for a relation between the new equations and a simpler linear equation of the type $ax - by = \pm 1$. This can be achieved either by raising the last equation to a certain power and equating the required variables with the corresponding expressions, or by using the method of hidden parameters. This method is as follows. Instead of each quantity included in the equation under study, we introduce an additional number of parameters using a linear or non-linear function, which allows us to find connections between the parameters, as well as some parametric solution of the desired equation. After that, we can reset (nullify) the introduced additional linear or non-linear dependence (in fact, hidden parameters), which in the original formulation of the problem would mean just replacing the letter designation of the parameters. However, as a result of these actions, a simplified parametric solution is obtained and the found relationship between the parameters remains. But if we do not reset (not nullify) the introduced additional parameters, then the solution turns out to be more cumbersome, with a large number of parameters, but including a larger number of parametric solutions. In any case, the solution is found through the solution of a linear equation of the type $ax - by = \pm 1$.

The proposed solution method was demonstrated in detail when solving the Pell equation, and also developed in a more general form and applied to the parametric solution of high-order Diophantine equations, such as the Delaunay equation, the Nagell equation and similar equations.

**Acknowledgments.** The author thanks Dr. Alexey Kuznetsov for careful reading, identifying typos, and providing valuable advice.

The main part of this material is an original translation of an article published in the journal "Trends in the Development of Science and Education." [in Russian]

## References

1. Arteha S.N.: Method of hidden parameters and Pell's equation. JP Journal of Algebra, Number Theory and Applications **2**(1), 21-46 (2002) http://doi.org/10.5281/zenodo.3463927
2. Andreescu T., Andrica D.: Quadratic Diophantine Equations. Springer, New York (2015) https://doi.org/10.1007/978-0-387-54109-9
3. Delaunay B.N., Faddeev D.K.: The theory of irrationalities of the third degree. Amer. Math. Soc. Providence, Rhode Island (1964) https://www.mathnet.ru/eng/tm/v11/p3
4. Jia Ch., Matsumoto K.: Analytic Number Theory. Springer, New York (2002) https://doi.org/10.1007/978-1-4757-3621-2
5. Matthews K.R., Robertson J.P., Srinivasan A.: On fundamental solutions of binary quadratic form equations. Acta Arithmetica **169**(3), 291-299 (2015) http://doi.org/10.4064/aa169-3-4
6. Nagell T.: Contributions to the theory of a category of Diophantine equations of the second degree with two unknowns. Nova Acta Regiae Soc. Sci. Upsaliensis **16**(2), Ser. 4, 38 (1955)
7. Nagell T.: Introduction to Number Theory. AMS Chelsea Publishing **163**, New York (1964) https://archive.org/details/introductiontonu0000nage/page/n5/mode/2up
8. Scherr Z., Thompson K.: Quartic integral polynomial Pell equations. Journal of Number Theory **259**, 38–56 (2024) https://doi.org/10.1016/j.jnt.2023.12.009
9. Togbé A., Walsh P.G.: A classical approach to a parametric family of simultaneous Pell equations with applications to a family of Thue equations. Bol. Soc. Mat. Mex. **28**, 4 (2022) https://doi.org/10.1007/s40590-021-00393-5
10. Tzanakis N., de Weger B.M.M.: On the practical solution of the Thue equation. J. Number Theor **31**(2), 99-132 (1989) https://doi.org/10.1016/0022-314X(89)90014-0